\documentclass[letterpaper,10pt,onecolumn]{cssconf}
\IEEEoverridecommandlockouts
\usepackage{amsfonts,amssymb}
\newcommand{\RR}{\rm I \kern -0.2em R}

\newtheorem{theorem}{Theorem}
\newtheorem{remark}{Remark}

\newtheorem{corollary}{Corollary}

\usepackage{epsfig}
\usepackage{cite}
\usepackage{color}
\usepackage{graphicx,color}
\usepackage[reqno]{amsmath}
\usepackage{amsbsy}
\usepackage{amssymb}

 \DeclareMathOperator{\sinc}{sinc}

\begin{document}
\thispagestyle{empty} \title{On the Stability of the Kuramoto
Model of Coupled Nonlinear Oscillators$^\dagger$}
\author{ Ali Jadbabaie$^*$, Nader Motee$^*$, and Mauricio
Barahona$^\ddagger$
\thanks{$^*$ A. Jadbabaie and N. Motee are with the Department of
Electrical and Systems Engineering and GRASP Laboratory,
University of Pennsylvania, Philadelphia, PA 19104. {\tt
 email: \{jadbabai, motee\}@grasp.upenn.edu}. A. Jadbabaie's research
is supported in part by ARO/MURI DAAD19-02-1-0383, ONR/YIP-542371
and NSF-ECS-0347285.} \noindent\thanks{$\ddagger$ M. Barahona is
with the Department of Bioengineering, Imperial College London,
United Kingdom.{\tt email: m.barahona@imperial.ac.uk}}
\thanks{$\dagger$ An earlier version of this paper was presented at
the 2004 American Control Conference, ACC 2004.}}
\date{}
\maketitle

\begin{abstract} We provide an analysis of the classic Kuramoto
model of coupled nonlinear oscillators that goes beyond the
existing results for all-to-all networks of identical oscillators.
Our work is applicable to oscillator networks of arbitrary
interconnection topology with uncertain natural frequencies. Using
tools from spectral graph theory and control theory, we prove that
for couplings above a critical value, the synchronized state is
locally asymptotically stable, resulting in convergence of all
phase differences to a constant value, both in the case of
identical natural frequencies as well as uncertain ones. We
further explain the behavior of the system as the number of
oscillators grows to infinity.
\end{abstract}

\section{Background and introduction}

Over the past decade, considerable attention has been devoted to
the problem of coordinated motion of multiple autonomous agents. A
variety of disciplines (as diverse as ecology, the social
sciences, statistical physics, computer graphics and, indeed,
systems and control theory) are developing an understanding of how
a group of moving objects (such as flocks of birds, schools of
fish, crowds of people~\cite{Vicnature,Natview}, or collections of
autonomous robots or unmanned
vehicles~\cite{CDC_boids_smooth,CDC_boids_nonsmooth}) can reach a
consensus and move in formation without centralized coordination.
Interestingly, this has coincided with a surge of activity in the
area of network dynamics, which focusses on the relationship
between graph structure and dynamical behavior of large networks
of diverse origin.

A classic example of distributed coordination in physics,
engineering and biology is the synchronization of arrays of
coupled nonlinear oscillators~\cite{strogatz_rev, strogatz94,
winfree1}. Building on long-standing experiments (dating back to
Huyghens and van der Pol), the problem of collective
synchronization was explored mathematically by the Russian school
of Andronov. Norbert Wiener~\cite{wiener1} also recognized its
ubiquity in the natural world, and even speculated about its
relevance to the existence of characteristic rhythms in the
brain~\cite{strogatz_book}.

Following on key insights by Winfree~\cite{winfree1},
Kuramoto~\cite{kuramoto1} proposed in the 1970s a tractable model
for oscillator synchronization that has become archetypal in the
physics and dynamical systems literatures.
(See~\cite{strogatz_rev} for an excellent review of the
state-of-the-art on this model.) More recently, researchers in the
control community~\cite{JusthKrisna02,Slotine,LucMureau} have
recognized that nonlinear synchronization phenomena are
mathematically related to the problem of coordination and
consensus among multi-agent systems~\cite{OlfMur,JaLiMo02}.

\section{Model description}
The classic Kuramoto model describes the dynamics of a set of $N$
phase oscillators $\theta_i$ with natural frequencies $\omega_i$.
The time evolution of the $i$-th oscillator is given by:
\begin{equation}
\dot \theta_i = \omega_i + \frac{K}{N} \sum_{j=1}^N \sin(\theta_j
- \theta_i), \label{eq:dynamics}
\end{equation}
where $K$ is the coupling strength, a key parameter in the
problem. One of Kuramoto's results was to show numerically that
when the $\omega_i$'s are randomly chosen from a Cauchy
probability distribution in the infinite $N$ limit, there is a
critical value of the coupling above which all phase differences
remain constant, i.e., the oscillators
synchronize~\cite{kuramoto1,kuramoto2}. If we think of the
oscillators as points moving on a circle, they would rotate
keeping the phase differences constant.

Kuramoto used the magnitude $R$ of the centroid of the points as a
`natural' measure of synchronization:
\begin{equation}
R \,e^\psi= \frac{1}{N}\sum_{i=1}^N e^{j \,\theta_i}.
\label{eq:orderparameter}
\end{equation}
Clearly, if all the $\omega_i$'s are the same then $R=1$ when all
agents are in sync. If the natural frequencies are not identical
but the oscillators synchronize, $R$ converges to a
\textit{constant} $R_\infty < 1$. On the other hand, when all
agents are completely out of phase with respect to each other the
value of $R$ remains close to $0$ most of the time. Because it
characterizes the dynamical behavior of the system, $R$ is
referred to as the {\it order parameter} in the physics
literature.

Kuramoto's analysis used simple trigonometry to rewrite the state
equation~(\ref{eq:dynamics}) in terms of the order parameter.
After switching to a rotating frame,
Eq.~(\ref{eq:dynamics}) becomes:
\begin{equation}
\dot \theta_i = \omega_i -\frac{K}{N} \, R \sin(\theta_i - \psi).
\label{eq:all-to-all-simple}
\end{equation}
In other words, each phase is modulated by the magnitude $R$ and
phase $\psi$ of the \textit{average} phasor.
In physics notation, this constitutes a {\it mean field} or
``all-to-all" model.

With some brilliant intuition, Kuramoto showed that for an
infinite number of oscillators there is a critical coupling $K_c$
below which the oscillators are incoherent (i.e., fully
unsynchronized). In addition, there is another critical coupling
$K_L \geq K_c$ above which \textit{all} oscillators are
synchronized. In that regime, the order parameter $R$ grows
exponentially in time until it saturates at a value $R_\infty(K)
\leq 1$. The branch of $R$ with $K>K_L$ is called the fully
synchronized state, while $K<K_c$ corresponds to the totally
unsynchronized state. Kuramoto also calculated analytically the
value for $K_c$ and $ R_\infty$ for a few well-known distributions
in the case of an infinite number of oscillators connected
all-to-all.

Despite its success, several aspects of the well-studied $N \to
\infty$, all-to-all Kuramoto model are still a puzzle, as
summarized beautifully in the review by Steve
Strogatz~\cite{strogatz_rev}. For instance, what does it mean that
$R$ stays close to zero in the unsynchronized state $K<K_c$? This
cannot be true at all times: when $K=0$ and the $\omega_i$'s are
irrational with respect to each other, the trajectories are dense
on the $N$-torus resulting in an $R$ which will almost surely
visit any number between 0 and 1. However, simulations indicate
that it is true most of the time. On the other extreme, the case
of few oscillators has been tackled in the dynamical systems
literature with rigorous bifurcation analysis. However, even basic
results are not available for the large but finite $N$ case, which
is of utmost interest in systems engineering.

Our goal here is to perform a system theoretic analysis of the
finite $N$ case with arbitrary connectivity. To proceed, we
rewrite the model in terms of the incidence matrix of the
undirected graph that describes the interconnection topology---
the standard all-to-all case is then the specific case of the
complete graph. We then provide several necessary as well as
sufficient lower bounds for the critical coupling $K_L$. These
include a bound for $K$ below which there is no fixed-point, and a
value of $K$ above which there is a unique fixed-point.
This extends similar results in~\cite{VanHemm93} for the case of 2
oscillators with a finite set of values for the natural
frequencies.

\section{Graph theoretical formulation of Kuramoto's model}

A good source for the necessary graph theory terminology
is~\cite{GoGo01}. We formalize our results through two matrices
that encode the topology of the connections. The incidence matrix
$B$ of an oriented graph $\mathcal{G}^\sigma$ with $N$ vertices
and $e$ edges is the $N \times e$ matrix such that: $B_{ij} = 1$
if the edge $j$ is incoming to vertex $i$, $B_{ij}=-1$ if edge $j$
is outcoming from vertex $i$, and $0$ otherwise. The symmetric $N
\times N$ matrix defined as:
$L = 
BB^T$
is called the Laplacian of $\mathcal{G}$ and is independent of the
choice of orientation $\sigma$. The Laplacian has several
important properties: $L$ is always positive semidefinite with a
zero eigenvalue; the algebraic multiplicity of its zero eigenvalue
is equal to the number of connected components in the graph; the
$N$-dimensional eigenvector associated with the zero eigenvalue is
the vector of ones, ${\bf 1}_N$. It is known that the spectrum of
the Laplacian matrix $\{\lambda_i (L)\}$ captures many topological
properties of the graph. Specifically, Fiedler showed that the
first non-zero eigenvalue $\lambda_2(L)$ (sometimes denoted the
algebraic connectivity) gives a measure of connectedness of the
graph. If we associate a positive number $W_i$ to each edge and we
form the diagonal matrix $W_{e \times e}:=\mbox{diag}(W_i)$, then
the matrix
$L_W(\mathcal{G}) = B W B^T$
is a weighted Laplacian which fulfills the above properties.

In this framework, the Kuramoto model~(\ref{eq:dynamics}) can be
generalized to any general interconnection topology as:
\begin{equation}
\label{eq:kuramoto_matrix}
  \dot \theta = \omega -\frac{K}{N} \, B \sin (B^T \theta),
\end{equation}
where $B$ is the incidence matrix of
the unweighted graph, and $\theta$ and $\omega$ are $N \times 1$
vectors. (It is also helpful to define the $e \times 1$ vector of
phase differences $\phi:= B^T \theta$.) A generalization of the
order parameter defined in~(\ref{eq:orderparameter}) for the
general Kuramoto model is:
\begin{equation}
\label{eq:orderparameter_matrix}
  r^2 = \frac{N^2 - 2 e+2 \, {\bf 1}_e^T \cos(B^T \theta)}{N^2}.
\end{equation}
It is easy to show that when the graph is complete, this is the
square of the magnitude of the average phasor, i.e., for $B=B_c$,
we have  $r_c^2= R^2$. While the average phasor interpretation
does not generalize to general connected topologies, the above
notion of an order parameter does generalize to arbitrary
connected graphs. In fact the order parameter can be conveniently
written in terms of the Laplacian matrix $L$ of the underlying
graph, as a measure of synchrony, or alignment. Specifically,
after some algebra, Equation~(\ref{eq:orderparameter_matrix}) can
be written as
\[
r^2 =1-\frac{1}{N} [e^{j \theta}]^* L [e^{j \theta}] =
1-\frac{1}{N}( [\cos \theta]^T L [\cos \theta] + [\sin \theta] ^T
L [\sin \theta]).
\]
where $[e^{j \theta}]= [e^{j \theta_1} \cdots e^{j \theta_N}]^T $
is the vector of complex phasors, and $^*$ denotes complex
conjugate transpose.

 The above equation provides us with an
interesting interpretation of the order parameter: each individual
oscillator $i$ can be thought of as a rotor moving on a circle
with unit radius, with velocity vector $v_i = e^{j \theta_i}$. The
total measure of disagreement between all velocity vectors, or the
global measure of asynchrony can be written as $\sum_{i} \sum_{ j
\in {\mathcal N}_i} ||v_i -v_j||^2$, which is nothing but
$N^2(1-r^2)$. In other words, the order parameter is a scaled
measure of {\it agreement} among velocity vectors of rotors, hence
a measure of synchronization. This allows us to extend the
stability analysis to a graph with arbitrary connected topology.

\begin{remark}
In the limit of small angles, the general Kuramoto
model~(\ref{eq:kuramoto_matrix}) gives the continuous-time Vicsek
flocking boid model~\cite{Vicsek1} which was analyzed
in~\cite{JaLiMo02}: $\dot \theta = \omega -(K/N) \, B \sin (B^T
\theta) \thickapprox
 \omega -(K/N)\, L \, \theta.$
Conversely, the classic Kuramoto model~(\ref{eq:dynamics}) can be
thought of as a nonlinear extension of the Vicsek model for a
complete graph.
\end{remark}

\begin{remark}
It is straightforward to show that the analytical simplification
(Eq.~\ref{eq:all-to-all-simple}) in the (standard) all-to-all
model appear as a result of the special symmetry of the Laplacian
of the complete graph:
\begin{equation}
L_c = B_c B_c^T = N I - {\bf 1}_N {\bf 1}_N^T.
\label{eq:Laplacian_complete}
\end{equation}
\end{remark}

\section{Synchronization of identical coupled oscillators}

We start by considering the general Kuramoto
model~(\ref{eq:kuramoto_matrix}) in its unperturbed version, i.e.,
when all the natural frequencies $\omega_i$ are identical:
\begin{equation}
\label{eq:identical}
 \dot \theta = - \frac{K}{N} \, B \sin(B^T \theta).
\end{equation}
(By switching to a rotating frame, it is easily shown that we can
assume that the natural frequencies $\omega_i$ are all zero,
without loss of generality.)
\begin{theorem}
\label{thm:kuram} Consider the unperturbed Kuramoto
model~(\ref{eq:identical}) defined over an arbitrary connected
graph with incidence matrix $B$. For  any value of the coupling $K
> 0$, all trajectories will converge to the set of equilibrium solutions. In particular the synchronized state is locally  asymptotically stable.
Moreover, the rate of approach to the synchronized state is no
worse than $ (2 \, K / \pi N) \, \lambda_2(L)$, where
$\lambda_2(L)$ is the Fiedler eigenvalue or the algebraic
connectivity of the graph.
\end{theorem}

\begin{proof}
Consider the function
$U_1(\theta)=
1 - r^2 = \frac{4 ||\sin(\frac{B^T \theta}{2})||^2}{N^2}$, where
$r^2$ has been defined in~(\ref{eq:orderparameter_matrix}).
A simple calculation reveals that $\nabla_\theta U = (2/N^2) B
\sin(B^T \theta)$ which leads to
\[ \dot U(\theta)=
\nabla_\theta U \,\, \dot\theta = - \frac{2}{KN} \dot \theta^T
\dot \theta \le 0.
\]
Therefore, the positive function $0 \leq U(\theta) \leq 1$ is a
non-increasing function along the trajectories of the system.  By
using LaSalle's invariance principle  we conclude that $U$ is a
Lyapunov function for the system, and that all trajectories
converge to the set where $\dot \theta$ is zero, i.e., the
equilibrium solutions.

Define now the $e \times e$ diagonal matrix $W(\phi) := \mbox{
diag}(\sinc(\phi_i))$, where $\sinc(\phi_i) = \sin(\phi_i)/
\phi_i$ is positive for $\phi_i \in (-\pi,\pi)^e$. Note also that
the angles move on an N- torus. Consider the proper subset of the
torus in which  $\phi \in (-\pi,\pi)^e$. The diagonal weight
matrix $W(\phi) > 0$ can be thought of as phase-dependent weight
functions on the graph. The trajectories converge to fixed-points,
which are the solutions of $ L_W \, \theta := \left (B \, W(\phi)
\, B^T \right) \theta = 0.$ The fact that $(\theta_0 {\bf 1}_N)$
is a locally exponentially stable equilibrium solution follows
easily: for any connected graph the null space of the weighted
Laplacian contains only  the vector ${\bf 1}_N$. Note that in
general, for an arbitrary connected topology the system has many
other equilibrium solutions, some of which might even be locally
stable. One such example is the  {\it ring} topology~\cite{naomi}.

Alternatively, one could use a different Lyapunov function similar
to the approach in~\cite{Olfati_acc03} and consider the quadratic
Lyapunov function candidate
$U_2 =\frac{1}{2} \sum_{i=1}^N\sum_{j=1}^N (\theta_i -\theta_j)^2
= \theta^T L_c \theta$, where $L_c = N I - {\bf 1 1^T}$  is the
Laplacian matrix of a complete graph. Note that $B^T{\bf 1} =0$,
therefore, a simple calculation reveals that
\begin{equation}
\dot{U_2}
= - \frac{K}{N} \theta^T B \sin (B^T \theta) = -\frac{K}{N} \theta
^T B W(\phi) B^T\theta \le 0. \nonumber
\end{equation}
It is interesting to note that the Lyapunov function $U_2$ is
small angle approximation  of Lyapunov function $U_1$.

Using the same argument as above, we conclude that the largest
sublevel set of $ U_2$ which is contained inside
$|\theta_l|<\frac{\pi}{2}\, l=1,\cdots,e$ is positively invariant.
 With the quadratic function $U_2(\theta)$
however, we can show that locally, the convergence is exponential
with the rate determined by the second smallest eigenvalue of the
weighted Laplacian:
\begin{equation}
\dot{U_2} \le -\frac{K}{N} \lambda_2(BW(\phi)B^T) ||\theta_{\bf 1^
\perp}||^2  \leq - \frac{2K}{\pi N} \lambda_2(L) ||\theta_{\bf 1^
\perp}||^2, \nonumber
\end{equation}
since $\lambda_2(BW(\phi)B^T) \leq (2/\pi) \lambda_2(BB^T)$.
\end{proof}

\begin{corollary}
For the complete graph, $\lambda_2(L_c) = N$ and the
synchronization rate for the mean-field model is no worse than $2
K / \pi$.
\end{corollary}

\begin{remark}
Similar results hold even if the topology of the graph changes in
time~\cite{JaLiMo02}. The result can be extended to general
notions of connectivity, i.e., when the interconnection graph is
not connected at all times but there is a path between any two
nodes over contiguous, non-overlapping, and uniformly bounded time
intervals. It is also possible to generalize to the case of
directed graphs by introducing notions of weak
connectivity~\cite{LucMureau}.
\end{remark}

\begin{remark}
The synchronization argument can be readily extended to the case
of more complicated coupling functions $f(\cdot)$ (other than the
$\sin(\cdot)$ function) so long as $\phi^T f(\phi) \ge 0$.
\end{remark}

\begin{remark}
The function ${\bf 1}_e^T \cos(B^T \theta)$
is an energy function for the XY-model in statistical physics. It
was considered as a Lyapunov-like function for the Kuramoto model
by Van Hemmen and Wreszinski~\cite{VanHemm93}, as well as
in~\cite{JusthKrisna02}.
\end{remark}

\begin{remark}
The global results obtained by Watanabe and
Strogatz~\cite{WatanabeStrogatz94} require all-to-all
connectivity. An extension of the methodology
in~\cite{WatanabeStrogatz94} to arbitrary topologies does not
appear to be trivial.
\end{remark}

\section{The case of non-identical oscillators}
In the rest of the paper we treat the more complicated case of
oscillators with non-identical natural frequencies. Although there
is an extensive literature for the $N \to \infty$ case with
all-to-all connectivity, we will focus here on the case of finite
$N$ and arbitrary topology given by
Eq.~(\ref{eq:kuramoto_matrix}). We consider the frequencies to be
random perturbations which, albeit drawn from a probability
distribution, remain constant in time, i.e., the
dynamics~(\ref{eq:kuramoto_matrix}) is deterministic yet
uncertain. This problem is distinct to some treatments in the
physics literature, which transform the problem into a
Fokker-Planck equation, effectively connected to a
\textit{stochastic} differential equation.


Synchronization is best defined in a \textit{grounded} system,
where the phases are defined with respect to a reference variable
(or 'ground').
This can be achieved by any projection $V_{N \times (N-1)}$ such
that
\begin{equation}\label{eq:projection}
  V^T V =I, \quad V V^T = I - \frac{ {\bf 1}_N {\bf
1}_N^T}{N} = \frac{L_c}{N}, \quad V^T{\bf 1}_N = {\bf 0}.
\end{equation}
Thus, $V$ is a matrix of $N-1$ orthonormal vectors orthogonal to
the vector ${\bf 1}_N$ which generate the set of grounded
coordinates $\bar \theta:=V^T \theta$ and frequencies $\bar
\omega:= V^T \omega$. The grounded Kuramoto model is:
\begin{equation}
\dot {\bar \theta}= \bar \omega -\frac{K}{N} V^T B \sin(B^T V \bar
\theta) = \bar \omega -\frac{K}{N} V^T B W(\bar \theta) B^T V \bar
\theta,  \label{eq:grounded}
\end{equation}
where, again, $W(\bar \theta):= \rm{diag}(\sinc(\phi_i))$ and
$\phi = B^T V \bar \theta$. In this grounded system, the
synchronized state is a \textit{fixed point}.

\begin{remark}
From Eq.~(\ref{eq:grounded}) it is easy to see why the natural
frequencies can be centered around zero without loss of
generality. Multiply Eq.~(\ref{eq:grounded}) from the left by $V$
and use~(\ref{eq:projection}) and $B^T {\bf 1}_N =0$ to recover
the original Eq.~(\ref{eq:kuramoto_matrix}) with new variables
$\Theta = \theta -  [\langle \omega \rangle t]\, {\bf 1}_N$ and
frequencies $\Omega = \omega - \langle \omega \rangle \, {\bf
1}_N$, where $<\omega>$ is the average frequency.
\end{remark}



\section{Bound for the asymptotic value of the order parameter}




Consider a Lyapunov function candidate based on the square of the
order parameter $r^2$ defined in~(\ref{eq:orderparameter_matrix}).
The derivative of this function along the trajectories is
\[
\dot{r^2} =\frac{1}{N^2}\left[\frac{K}{N} (\sin B^T \theta)^T B^T
B (\sin B^T \theta)- \omega^T B \sin B^T \theta \right],
\]
which is an ellipsoid in the $\sin (B^T \theta)$ coordinate
centered at $\frac{N \omega}{K}$. Outside of  a neighborhood of
the origin given by
\begin{equation}
\label{deriv_Lyap}
 ||B \sin(B^T \theta)||_2> \frac{N}{K}
||\omega||_2
\end{equation}
the derivative is positive, resulting in growth of the order
parameter. The boundary of this region contains the equilibria. By
using an ultimate boundedness argument, the trajectories are
confined to the smallest sublevel-set of $r$ containing the set
defined by~(\ref{deriv_Lyap}).

We now use~(\ref{deriv_Lyap}) to obtain an estimate of the
asymptotic value of the order parameter. The vector $\sin(B^T
\theta)$ can be decomposed into two orthogonal components:
$y_1(\theta)$, in the null space of $B$, and $y_2(\theta)$ in the
range space of $B^T$. The first component is annihilated when it
is multiplied by $B$ on the left. As  a result, the region over
which $\dot r^2$ is positive can be characterized as
\[
||y_2(\theta)||_2>\frac{N}{K\sqrt{\lambda_2(L)}} ||\omega||_2.
\]
where $\lambda_2(L)$ is the algebraic connectivity of the
unweighted graph. We now bound the value of $U$ over the region
where $\dot r^2$ is negative. A simple bounding reveals that
{\small
\[
2{\bf 1}^T \cos(B^T \theta) \le ||{\bf 1}||^2 + ||\cos(B^T
\theta)||^2 = 2||{\bf 1}||^2 -||\sin(B^T \theta)||^2,
\]}
from which {\small
\[
r^2 \le  \frac{N^2 -||\sin B^T \theta||^2}{N^2} \le \frac{N^2
-||y_2( \theta)||^2}{N^2} \le \frac{N^2 - \frac{N^2 ||w||^2}{K^2
\lambda_2(L)}}{N^2} .
\]}

We can immediately observe that the asymptotic behavior of the
order parameter is inversely proportional to the algebraic
connectivity of the graph. Of course, because of the
over-bounding, the bound is conservative---its asymptotic value is
1 as opposed to the actual less-than-one value. Nevertheless, this
gives us a bound on the growth rate of $r^2$, and, as a result,
the growth rate on $r$ is bounded by
$\frac{1}{\sqrt{\lambda_2(L)}}$.

This means that asymptotically
\[
r \le \sqrt{ 1 - \frac{||w||^2}{K^2 \lambda_2(L)}}
\]
which would result in an increase rate of
$\mathcal{O}(\frac{1}{\sqrt{N}})$ when the graph is {\it
complete}.
\begin{remark}
Consider a complete graph where the natural frequencies are
independent random variables chosen from a normal distribution
$\omega_i \; \sim \mathcal{N}(0,\sigma)$.  Then $||w||_2$ scales
as $\sqrt{N} \sigma$, which results in a bound for $r <
\sqrt{1-(\sigma/K)^2}$ that is independent of $N$.
\end{remark}

\begin{remark}
In~\cite{VanHemm93}, the authors added a linear term $\omega^T
\theta$ to the Lyapunov function candidate to guarantee negativity
of the derivative everywhere except at the fixed-points, reducing
the perturbed model to a gradient system. The linear term,
however, makes the Lyapunov function indefinite.
\end{remark}

We will see in the next section that if $K$ is large enough to
guarantee the existence of a unique fixed point (via a contraction
argument), condition~(\ref{deriv_Lyap}) will be trivially
satisfied. This means that if $K$ is large enough the derivative
of the order parameter will be positive, resulting in the
asymptotic stability of the synchronized state.

\section{Bounds for the critical coupling}

As the coupling $K$ is decreased, there is a critical value $K_L$
below which no fixed point exists, resulting in a running solution
for the grounded system~(\ref{eq:grounded}). This means that the
system cannot be fully synchronized for $K < K_L$.

An easy \textit{sufficient} condition for the fixed point $\bar
\theta^*$ to be stable is  for $\phi^* = B^T V \bar \theta^*$ to
be contained  in any closed subset of
$(-\frac{\pi}{2},\frac{\pi}{2})^{e}$, which implies that
$|\theta^*| < \frac{\pi}{4}$.  This is demonstrated by taking the
Jacobian of $V^T B \sin B^T \theta$, and noting that it is equal
to $V^T B\mbox{diag}[\cos(B^T \theta^*)]B^T V$, which is positive
definite over that set.

\subsection{Critical value of coupling for complete graphs}
Our results generalize those of Van Hemmen {\it et
al.}~\cite{VanHemm93} in the case of a complete graph.
Specifically, it can be shown that the critical value of the
coupling is determined by the value of $K$ for which  the fixed
point disappears. This can be explained by looking at the fixed
point equation $ B \sin(B^T  \theta^\ast) = \frac{N \omega}{K}.$

Let $\omega_{max}=||\omega||_{\infty}$ and note that the induced
infinity norm of a matrix is the maximum absolute row sum, i.e.,
$||B||_\infty = d_{max}$, where $d_{max}$ is the maximum degree of
the graph. In the case of a complete graph, $d_{max}= N-1$. Then,
\[
\frac{N \omega_{max}}{K} \le d_{max}
\]
resulting in the following lower bound for $K_L$, the coupling
above which a fixed point exists:
\[
K_L > \frac{N \omega_{max}}{d_{max}}.
\]

This bound can be tightened by using the generalized inverse of
$V^T B$ and bounding the component of the $\sin (B^T \theta)$ in
the range of $B^T$. The generalized inverse, denoted by $(V^T
B)^{\#}$, is equal to $B^T V \Lambda^{-1}$, where $\Lambda$ is the
$N-1$ diagonal matrix of the eigenvalues of the unweighted
Laplacian. We therefore have the following expression
\[(\sin(B^T\theta))_{R(B^T)} = B^T V \Lambda^{-1}V^T  \frac{N
\omega}{K}.
\]
Noting that $ L^\# =V \Lambda^{-1}V^T$, we have
\[
(\sin(B^T\theta))_{R(B^T)} = B^T L^\# B \sin(B^T \theta)= B^T L^\#
\frac{N \omega}{K}.
\]

The generalized inverse of the Laplacian, in the case of a
complete graph can be written as $L_c^\# = \frac{1}{N}(I -
\frac{\bf 1 1 ^T}{N})$. Noting that the infinity norm of the
$\sin$ vector is less than or equal to 1, and that $B^T L^\# B =
\frac{B^T B}{N}$, we have
\[
\frac{||B^T \omega||_{\infty}}{K} \le \frac{||B^T
B||_{\infty}}{N},
\]
which gives us the bound
\[
K_L \ge ||B^T \omega||_{\infty} \frac{N}{2(N-1)}.
\]

This is in excellent agreement with that of Van Hemmen {\it et
al.}~\cite{VanHemm93} which they obtained for the simplest case of
two oscillators.

\begin{remark}
If the graph is a tree, $V^T B$ has full row rank and $\sin(B
\theta)$ does not have a component in the null space of $L$. In
that case $K_L > ||B^T L^\# \omega||_{\infty}$ is a tight bound,
meaning that it is necessary and sufficient for synchronization.
In the general case, however, this bound is just necessary.
\end{remark}

\subsection{Existence and uniqueness of stable fixed points}

The fixed point equation can be written as
\[
  \theta = ( B W(B^T \theta)B^T )^{\#} \frac {N
\omega}{K} = L_W^\#(B^T \theta) \frac{N \omega}{K}.
\]
Using Brouwer's fixed point theorem (i.e., a continuous function
that maps a non-empty compact, convex set $X$ into itself has at
least one fixed-point), we can develop conditions which guarantee
the existence (but not uniqueness) of the fixed point. If a
fixed-point exists in any compact subset of $ \theta \in
(-\frac{\pi}{4}, \frac{\pi}{4})$, it is stable, since this will
ensure that $B^T \theta$ is between $-\frac{\pi}{2}$ and
$\frac{\pi}{2}$. We therefore have to ensure that
\[
K > \frac{4}{\pi}N \max_{|\theta_i|<\frac{\pi}{4}} ||L_W^\# (B^T
\theta)||_{\infty} ||\omega||_\infty.
\]
Simulations indicate that in the case  of a complete graph, the
infinity norm of the matrix $L_W^\#$ scales as
$\mathcal{O}(\frac{1}{N})$. It is worth mentioning that the norm
of $L_W^\#$ is a well studied object in the theory of Markov
chains. The infinity norm of $L_W^\#$ is  a measure of the
sensitivity of the stationary distribution of the chain associated
with $L$ with respect to additive perturbations~\cite{ChoMey01}.

If the uncertain natural frequencies are 2-norm bounded, a better
strategy would be to impose the boundedness condition with respect
to the Euclidean norm. A sufficient condition for local stability
of the fixed-point is for $\theta_i$ to belong to
$(-\frac{\pi}{4}, \frac{\pi}{4})$. This amounts to having the
Euclidean norm of $\theta$ be less than $\frac{\pi}{4}\sqrt{N}$.
Again, using Brouwer's sufficient condition for existence of
fixed-points we have:

\[
 ||(B W(B^T \theta)B^T )^{\#}||_2 \frac {N
||\omega||_2}{K} \le \frac{\pi}{4}\sqrt{N}.
\]
Hence, a sufficient condition for synchronization of all
oscillators can be determined in terms of a lower bound for $K$:
\[
K_L \ge \frac{4}{\pi} \frac{\sqrt{N} ||w||_2}{\min_{|\theta_i| \le
\frac{\pi}{4}} \lambda_2 (L_W(\theta))},
\]
where we used the fact that $||(BW(B^T
\theta)B)^\#||_2=\frac{1}{\lambda_2(L_W)}$, and $\lambda_2$ is the
algebraic connectivity of the (weighted) graph. A lower bound on
the minimum value of $\lambda_2$ occurs for the minimum value of
the weight which is $\frac{2}{\pi}$. As a result,
\begin{equation}
\label{eq:Kc} K_L \ge  2 \frac{\sqrt{N} ||w||_2}{\lambda_2(L)}.
\end{equation}

\begin{remark}
Using the upper bound provided for the order parameter earlier, we
can derive an upper bound for the asymptotic value of $r$ at
$K_L$: $r_\infty(K_L) \leq \frac{\sqrt{3}}{2}$. Furthermore, if
the stable fixed-point is in $(-\pi/4,\pi/4)^N$, then the order
parameter is lower bounded by $\sqrt{16-\pi^2}/4$. 
\end{remark}

\subsection{Bounds for the existence of a unique fixed-point}
In order to guarantee the existence of a unique fixed point we use
Banach's contraction principle and ensure that the right hand side
is a contraction. By noting that the Lipschitz constant for the
$\sinc(\cdot)$ function is $\alpha_s = \frac{1}{2}$,  we provide a
sufficient condition for contractivity (and therefore uniqueness
of the fixed-point).

We impose the contractivity condition on the $N-1$ dimensional
grounded system. In the grounded case, we have $\bar \theta =V^T
\theta$, and
\[
\bar \theta = (V^T  B W(B^T \theta)B^T V)^{-1} \frac {N V^T
\omega}{K}.
\]

After some algebra, the contraction requirement amounts to
\begin{equation}
\label{eq:Kc_cont} K \ge \frac{\pi^2}{4} \frac{N \lambda_{max}(L)
||w||_2}{\lambda_2(L)^2},
\end{equation}
where $\lambda_{max}$ is the largest eigenvalue of the Laplacian
of the graph.

Interestingly, this value of $K$ also ensures that the derivative
of $r^2$ is increasing,i.e., inequality~(\ref{deriv_Lyap}) is
satisfied, which means  that the order parameter is increasing. Of
course this is probably stronger than what is necessary for
uniqueness, as the contraction argument is only sufficient.
Nevertheless, we see that there is a large enough but finite value
of the coupling which guarantees the existence and uniqueness of
fixed points.

We now state the following theorem:
\begin{theorem}
Consider the Kuramoto model for non-identical coupled oscillators
with different natural frequencies $\omega_i$. For $K \ge K_L := 2
\frac{\sqrt{N} ||w||_2}{\lambda_2(L)}$, there exist at least one
fixed-point for $|\theta_i| <\frac{\pi}{4}$ or $|(B^T \theta)_i|
<\frac{\pi}{2}$. Moreover, for $K \ge \frac{\pi^2}{4} \frac{N
\lambda_{max}(L) ||w||_2}{\lambda_2(L)^2}$ there is only one
stable fixed-point (modulo a vector in the span of ${\bf 1}_N$),
and the order parameter is strictly increasing.
\end{theorem}
\begin{proof}
Proof of Theorem 3 : As we see before fixed point equation can be
written as
\begin{equation}
V^T B W(\bar \theta) B^T V \bar \theta = \frac{K}{N} \bar \omega
\end{equation}
or
\begin{eqnarray}
\bar \theta = (V^T B W(\bar \theta) B^T V)^{-1} \frac{K}{N} \bar \omega \\
=L_W (\bar \theta)^{-1} \frac{K}{N} \bar \omega \label{fixedpnt}
\end{eqnarray}

We will use Banach's contraction principle to show that
(\ref{fixedpnt}) has a unique solution when $V \bar \theta$ is in
any compact subset of $(-\frac{\pi}{2},\frac{\pi}{2})^{N}$. We
therefore need to show that
\begin{equation}
||\frac{N}{K} (L_W (\bar \theta_1)^{-1}-L_W (\bar
\theta_2)^{-1})\bar\omega||\leq \alpha ||\bar \theta_1-\bar
\theta_2||
\end{equation}
holds for some $0 \leq \alpha < 1$ and some norm. Using the
2-norm, we have
\begin{eqnarray}
||\frac{N}{K} (L_W (\bar \theta_1)^{-1}-L_W (\bar
\theta_2)^{-1})\bar\omega||_{2} &=& ||\frac{N}{K} L_W (\bar
\theta_1)^{-1}(L_W (\bar \theta_2)-L_W (\bar \theta_1))L_W
(\bar \theta_2)^{-1}\bar\omega||_{2}\nonumber\\
&\leq& \frac{N}{K}||L_W (\bar \theta_1)^{-1}||_{2} ||L_W (\bar
\theta_2)^{-1}||_{2} ||L_W (\bar \theta_2)-
L_W (\bar \theta_1)||_{2}||\bar\omega||_{2}\nonumber \\
&\leq& \frac{N}{K} \frac{1}{\lambda_{min} (L_W (\bar
\theta_1))}\frac{1}{\lambda_{min}(L_W (\bar \theta_2))}||V^T B
 (W(\bar\theta_1)-W(\bar\theta_2)) B^T V||_2 ||\omega||_2 \nonumber\\
&\leq& \frac{N}{K} \frac{1}{\lambda_{min}(L_W (\bar \theta_1))}\frac{1}{\lambda_{min}(L_W (\bar \theta_2))}||V^T B||_2 ||W(\bar\theta_1)-W(\bar\theta_2)||_\infty ||B^T V||_2 ||\omega||_2 \nonumber\\
&\leq& \frac{N}{K} \frac{1}{\lambda_{min}(L_W (\bar \theta_1))}\frac{1}{\lambda_{min}(L_W (\bar \theta_2))} \lambda_{\max}(L) ||W(\bar\theta_1)-W(\bar\theta_2)||_\infty ||\omega||_2 \nonumber\\
&\leq& \frac{N}{K} \frac{1}{\lambda_{min}(L_W (\bar
\theta_1))}\frac{1}{\lambda_{min}(L_W (\bar \theta_2))}
\lambda_{\max}(L) \alpha_{s} ||B^T||_{\infty}
||\bar\theta_1-\bar\theta_2||_{\infty} ||\omega||_2
\label{ineq_bound}
\end{eqnarray}

$\alpha_{s}$ is the Lipschitz constant of $sinc(.)$ which is
$0.5$. For a weighted graph with non-negative weights
$W=[w_{ij}]$, the second smallest non-zero eigenvalue is given by
~\cite{FChung}

\begin{equation}
\lambda_{min}(L_W (\bar \theta))=\lambda_{2}(L_W (\theta))=(N-1)
\inf_{f\perp D \bf 1} \frac{\sum_{i \sim j} (f_{i}-f_{j})^2
w_{ij}}{\sum_{v} f_{i}^2 d_{i}}
\end{equation}

with $d_{i}=\sum_{j} w_{ij}$ and $D$ denotes the diagonal matrix
with the $(i,i)$-th entry having value $d_{i}$. Since $B^T \theta
\in (-\frac{\pi}{2} , \frac{\pi}{2})^e$ ,we have $\frac{2}{\pi}
\leq  w_{ij} \leq 1 $. Therefore,

\begin{equation}
\lambda_{2}(L_W (\theta)) \geq \frac{2}{\pi} \lambda_{2}(L)
\label{egnvalue_bound}
\end{equation}
Applying (\ref{egnvalue_bound}) to (\ref{ineq_bound}), we get
\begin{eqnarray}
||\frac{N}{K} (L_W (\bar \theta_1)^{-1}-L_W (\bar
\theta_2)^{-1})\bar\omega||_{2} \leq \frac{N}{K} (\frac{\pi}{2})^2
\frac{\lambda_{N}(L)}{\lambda_2(L)^2}
||\bar\theta_1-\bar\theta_2||_{\infty} ||\omega||_2
\end{eqnarray}
 From here, we can find value of $K$'s that make this mapping contractive as follows
\begin{equation}
\frac{N}{K} (\frac{\pi}{2})^2
\frac{\lambda_{N}(L)}{\lambda_2(L)^2} ||\omega||_2 < 1
\end{equation}
or
\begin{equation}
 K > (\frac{\pi}{2})^2 \frac{N \lambda_{N}(L)}{\lambda_2(L)^2} ||\omega||_2
\end{equation}

As a result for $K>K_c$, where $K_c=(\frac{\pi}{2})^2 \frac{N
\lambda_{N}(L)}{\lambda_2(L)^2} ||\omega||_2$, the fixed-point
equation (\ref{fixedpnt}) has a unique and stable solution.
{\hfill \large{$\Box$}}
\end{proof}
%
%
\section{Concluding remarks}
In this paper we provided a stability analysis for the Kuramoto
model of coupled nonlinear oscillators for arbitrary topology. We
showed that when the oscillators are identical, there are at least
two Lyapunov functions which prove asymptotic stability of the
synchronized state, when all the phase differences are bounded by
$\frac{\pi}{2}$. We also showed that when the natural frequencies
are not the same, there is a critical value of the coupling below
which a totally synchronized state does not exist. Several bounds
for this critical value based on norm bounded uncertain natural
frequencies were shown to be in excellent agreement with existing
bounds in the physics literature for the case of the all-to-all
graph.

We also point out that contrary to the infinite $N$ case, there is
no partially synchronized state, i.e., for values of the coupling
below the critical value, the system of differential equations has
a running solution. Furthermore, we showed that there is always a
large enough but finite value of the coupling which results in
synchronization of oscillators and convergence of the angles to a
unique fixed-point. Another result of this paper is that the value
of the order parameter is not zero for the critical coupling
$K_L$. In fact, at least when the fixed-point is in the $(-\pi/2,
\pi/2)$ region, a rough estimate indicates that the value of $r$
is bounded between $\frac{\sqrt{16 - \pi^2}}{4}\approx 0.62$ and
$\frac{\sqrt{3}}{2}.$ Future research in this direction is needed
to determine the bound for $K$ when the natural frequencies are
not just norm bounded quantities but uncertain numbers chosen from
a probability distribution. Finally we mention that our value for
the upper bound of the order parameter is actually quite close to
simulations.

Our work hints at the advantageous marriage of systems and control
theory and graph theory, when studying dynamical systems over or
on networks~\cite{syncSW}.

\section{Acknowledgements}
We would like to thank Steve Strogatz, Naomi Leonard, George
Pappas and Bert Tanner for their helpful comments.

\bibliographystyle{plain}
\bibliography{kuramotobib_reallyfinal}


\end{document}